# The Age or Excess of the $M|G|\infty$ Queue Busy Cycle Mean Value

Manuel Alberto M. Ferreira[1], Marina Andrade[1] & José António Filipe[1]

[1] Instituto Universitário de Lisboa (ISCTE - IUL), BRU - IUL, Lisboa, Portugal

Correspondence: Manuel Alberto M. Ferreira, Instituto Universitário de Lisboa (ISCTE - IUL), BRU - IUL, Lisboa, Portugal. E-mail: manuel.ferreira@iscte.pt

**Abstract**

In this work, approximate and exact results concerning a performance measure of the $M|G|\infty$ system, the age or the excess of the busy cycle, are presented. It will be seen that it is also a measure of the busy period performance. Service distributions for which it is given for a simple expression and others for which this does not happen are considered. For this last case, bounds are deduced. A special emphasis is given to the exponential distribution and to those related with it, useful in reliability theory.

## 1. Introduction

In a queue system it is usual to call busy period a period that begins when a costumer arrives at the system being it empty, ends when a costumer abandons the system letting it empty, and during it there is always at least one customer being served.

So, in a queue system, there is a sequence of idle periods and busy periods.

Be then the $M|G|\infty$ system initially empty. The instants $0, t_1, t_2, \ldots$ at which the system enters in state 0, are a renewal process arrival instant, see (Hokstad, 1979). A cycle is complete whenever a renewal occurs, that is: an entrance at state 0. These cycles are called busy cycles and their length is a random variable designated Z.

In (Takács, 1962) it is showed that

$$E[Z] = \frac{e^\rho}{\lambda} \text{ and } E[Z^2] = \frac{2e^\rho \int_0^\infty \left(e^{-\lambda \int_0^t [1-G(v)]dv} - e^{-\rho}\right)dt}{\lambda} + \frac{2e^\rho}{\lambda^2} \quad (1)$$

where $\lambda$ is the Poisson process arrivals rate, $G(.)$ is the service distribution function, $\alpha$ its service time and $\rho = \lambda\alpha$ the traffic intensity.

Consider now a renewal process which time length between consecutive arrivals is a random variable Z and be $A(t)$ the time spent since the last renewal till *t*, or the time spent after *t* till the next renewal. If the renewals represent old devices turning out of order and being replaced, $A(t)$ is the *age* of the device in use at instant *t* or the remaining lifetime of a device in use at instant *t – excess* of a device in use at instant *t*, respectively.

Being interested in that device age, or excess, mean value, that is $\lim_{s\to\infty} \frac{\int_0^s A(t)dt}{s}$, it can be computed through, see (Ross, 1983).

$$\lim_{s\to\infty} \frac{\int_0^s A(t)dt}{s} = \frac{E[Z^2]}{2E[Z]} \quad (2)$$

Note that being the $M|G|\infty$ queue a system with no waiting and no losses, it is mandatory to present immediately an available server when every costumer arrives, it will be interesting, for instance, in a given instant of a busy period to have an idea of how much more time it will last. So there will have the notion of for how much time it will be necessary to have the servers available. This time is precisely the busy cycle excess.

The results presented will allow answering to this question in mean value terms.

## 2. $M|G|\infty$ Queue Busy Cycle Age or Excess Mean Value

Calling the $M|G|\infty$ system busy cycle age, or excess, mean value $\beta_c$.

$$\beta_c = \beta + \frac{1}{\lambda} \quad (3)$$

where

$$\beta = \int_0^\infty \left(e^{\rho - \lambda \int_0^t [1-G(v)]dv} - 1\right)dt \quad (4)$$

In (Sathe, 1985) it is showed that $\beta = \frac{1}{2\lambda}\rho^2(\gamma_s^2 + 1)E[e^{\lambda U(t)}]$ and $1 + 2\rho^{-2}(1 + \gamma_s^2)^{-1}\left(e^\rho - 1 - \rho - \frac{\rho^2}{2}\right) \leq E[e^{\lambda U(t)}] \leq 2\rho^{-2}(e^\rho - 1 - \rho)$ where $\gamma_s$ is the service distribution coefficient variation. So

$$\frac{\gamma_s^2}{\lambda}\frac{\rho^2}{2} - \alpha \leq \beta_c - E[Z] \leq \frac{\gamma_s^2}{\lambda}\sum_{n=2}^{\infty}\frac{\rho^n}{n!} - \alpha \tag{5}$$

The Expressions (3) and (4) show that $\beta_c$ depends on the whole service distribution structure and so is highly sensible to its form. The bounds given by (5) possess the great advantage of being valid whichever the service distribution is and depend only on $\rho, \lambda$ and $\gamma_s$.

The next proposition, immediate consequence of (5), allows to compare $\beta_c$ with $E[Z]$, that is insensible to the service distribution, since $\rho$ and $\gamma_s$ are known.

**Proposition 1**

$$If\ \gamma_s^2 \leq \frac{\rho}{e^\rho - 1 - \rho} \quad \beta_c \leq E[Z]$$
$$If\ \gamma_s^2 \geq \frac{2}{\rho} \quad \beta_c \geq E[Z]$$

**Observation:**

$$\frac{2}{\rho} \geq \frac{\rho}{e^\rho - 1 - \rho}$$

But $\lim_{\rho \to 0} \frac{2}{\rho} - \frac{\rho}{e^\rho - 1 - \rho} = \frac{2}{3}$.

## 3. Values of $\beta_c$ for Some Service Distributions

As it was emphasised in section 2, $\beta_c$ depends on the whole service time distribution. Then the values of $\beta_c$ for some service time distributions, obtained after (3), are presented:

**Exponential**

$$\beta_c^M = \frac{1}{\lambda} + \alpha \sum_{n=1}^{\infty} \frac{\rho^n}{nn!}$$

**Constant**

$$\beta_c^D = E[Z] - \alpha$$

$G(t) = \frac{e^{-\rho}}{e^{-\rho} + (1 - e^{-\rho})e^{-\lambda t}}, t \geq 0$, see (Ferreira & Andrade, 2012)

$$\beta_c = E[Z]$$

$G(t) = 1 - \frac{1}{1 + e^{-\rho}\left(e^{\frac{\lambda}{1-e^{-\rho}}t} - 1\right)}, t \geq 0$, see (Ferreira & Andrade, 2012)

$$\beta_c = \frac{e^\rho + e^{-\rho} - 1}{\lambda}$$

**Power function with parameter $c$ ($\alpha = \frac{c}{c+1}$)**

$$\beta_c^P = \frac{1}{\lambda} + \sum_{n=1}^{\infty}\sum_{k=0}^{n}\frac{\rho^{n-k}}{(n-k)!}\sum_{j=0}^{k}\frac{(-1)^j}{(c+1)^j j!(k-j)!(k+jc+1)}$$

**Uniform in $[0, 1]$ ($c=1$)**

$$\beta_c^P = \frac{1}{\lambda} + \sum_{n=1}^{\infty}\sum_{k=0}^{n}\frac{\rho^{n-k}}{2^{n-k}(n-k)!}\sum_{j=0}^{k}\frac{(-1)^j}{2^j j!(k-j)!(k+j+1)}$$

But note that the result for the constant service distribution may be derived easily from (5) making $\gamma_s = 0$.

For any service time distribution, after (5),

$$\beta_c \geq E[Z] - \alpha + \frac{\rho^2}{2\lambda}\gamma_s^2 \tag{6}$$

So, for the $M|G|\infty$ system, fixed $\alpha$ and $\lambda$, the least value of $\beta_c$ happens in the case of constant service time.

The values of $\beta_c$ are then computed for some values of $\alpha$ and $\lambda$ and presented in Tables 1 and 2.

Table 1. $\beta_c$ values for $\lambda = 1$ and various values of $\alpha$, for some service time distributions

| Service Time Distribution | $\alpha = .5$ | $\alpha = 1$ | $\alpha = 5$ | $\alpha = 10$ | $\alpha = 50$ |
|---|---|---|---|---|---|
| Exponential | 1.2850757 | 2.3178568 | 186.93907 | 24755.984 | $5.2920661 \times 10^{21}$ |
| Constant | 1.1487213 | 1.7182818 | 143.41316 | 22016.466 | $5.1847055 \times 10^{21}$ |
| $G(t) = \dfrac{e^{-\rho}}{e^{-\rho} + (1 - e^{-\rho})e^{-\lambda t}}$, $t \geq 0$ | 1.6487213 | 2.7182818 | 148.41316 | 22026.466 | $5.1847055 \times 10^{21}$ |
| $G(t) = 1 - \dfrac{1}{1 + e^{-\rho}\left(e^{\frac{\lambda}{1-e^{-\rho}}t} - 1\right)}$, $t \geq 0$ | 1.2552519 | 2.0861613 | 147.41990 | 22025.466 | $5.1847055 \times 10^{21}$ |

Table 2. $\beta_c$ values for $\alpha = 0.5$ and various values of $\lambda$, for some service time distributions

| Service Time Distribution | $\lambda = 2$ | $\lambda = 10$ | $\lambda = 20$ | $\lambda = 100$ |
|---|---|---|---|---|
| Exponential | 1.1589511 | 19.099311 | 1244.7304 | $5.9392749 \times 10^{19}$ |
| Constant | .85914091 | 14.341316 | 1100.8233 | $5.1847055 \times 10^{19}$ |
| $G(t) = \dfrac{e^{-\rho}}{e^{-\rho} + (1 - e^{-\rho})e^{-\lambda t}}$, $t \geq 0$ | 1.3591409 | 14.841316 | 1101.3233 | $5.1847055 \times 10^{19}$ |
| $G(t) = 1 - \dfrac{1}{1 + e^{-\rho}\left(e^{\frac{\lambda}{1-e^{-\rho}}t} - 1\right)}$, $t \geq 0$ | 1.0430806 | 14.741990 | 1101.2733 | $5.1847055 \times 10^{19}$ |
| Power function with parameter $c$ ($\alpha = \frac{c}{c+1}$) | 1.9626517 | 17.272158 | 168.2805 | $5.2381918 \times 10^{19}$ |

For the service time distributions exponential and power function the $\beta_c$ values were obtained through (3) by numerical methods.

The values in Tables 1 and 2 evidence the dependence of $\beta_c$ from the service time distribution structure, although for high values of $\rho$ that dependence vanishes expressively.

After (5) and the expression showed in this section for $\beta_c^M$ it is possible to obtain lower and upper bounds for this parameter, since the busy period mean value is $E[B] = \dfrac{e^{\rho}-1}{\lambda}$ for any service time distribution.

If the service time distribution is *NBUE*-**N**ew **B**etter than **U**sed in **E**xpectation with mean $\alpha$, $\int_t^\infty [1 - G(v)]dv \leq \int_t^\infty e^{-\frac{v}{\alpha}}dv$, see (Ross, 1983), and the lower bound obtained for $\beta_c^M$ is good for $\beta_c^{NBUE}$.

If the service time distribution is *NWUE*-**N**ew **W**orse than **U**sed in **E**xpectation with mean $\alpha$, $\int_t^\infty [1 - G(v)]dv \geq \int_t^\infty e^{-\frac{v}{\alpha}}dv$, see (Ross, 1983), and the upper bound obtained for $\beta_c^M$ is good for $\beta_c^{NWUE}$.

If the service time distribution is DFR-**D**ecreasing **F**ailure **R**ate, $1 - G(t) \geq e^{-\frac{t}{\alpha} - \frac{\gamma_s^2}{2} - \frac{1}{2}}$, see (Ross, 1983), and a lower bound for $\beta_c^{DFR}$ may be obtained following a methodology analogous to the one that allowed to obtain the lower bound for $\beta_c^M$. If the service time distribution is IMRL-**I**ncreasing **M**ean **R**esidual **L**ife, $1 - \dfrac{\int_0^t [1-G(v)]dv}{\alpha} \geq e^{-\frac{2t\alpha}{\mu_2} - \frac{2\alpha}{3\mu_2^2}\mu_3 + 1}$, being $\mu_r$ the $G(.)$ $r^{\text{th}}$ order moment around the origin, see (Brown, 1981) and (Cox, 1962), and it is possible to find a lower bound for $\beta_c^{IMRL}$ analogous to the one for $\beta_c^M$. For the power function service distribution, as $\gamma_s^2 = [c(c+2)]^{-1}$, a lower bound and an upper bound for $\beta_c^P$ that, for $c = 1$, are also valid for the uniform in $[0,1]$ service distribution are easily obtained. So

- **$\beta_c$ lower bounds**

**a) M and NWUE**

$$E[Z] - \alpha + \frac{\alpha\rho}{2}\left(1 + \frac{\rho}{6}\right)$$

**b) DFR**

$$E[Z] - \alpha + \frac{\alpha\rho}{2}\left(2e^{\frac{1+\gamma_s^2}{2}} - 1 + \rho\frac{3e^{1-\gamma_s^2} - 2}{6}\right)$$

**c) IMRL**

$$E[Z] - \alpha + \frac{\lambda}{4}\left(\mu_2 e^{1-\frac{2\alpha}{3\mu_2^2}\mu_3} - 2\alpha^2 + \rho\frac{\mu_2 e^{2(1-\frac{2\alpha}{3\mu_2^2}\mu_3)} - 4\alpha^2}{6}\right)$$

**d) Power function with parameter $c$**

$$E[Z] + \frac{\rho - 2c(c+2)}{2(c+1)(c+2)}$$

**e) Power function with parameter $c$**

$$E[Z] + \frac{\rho - 2c(c+2)}{2(c+1)(c+2)}$$

**f) Uniform in $[0, 1]$ $(c = 1)$**

$$E[Z] + \frac{\frac{\lambda}{2} - 6}{12}$$

- **$\beta_c$ upper bounds**

  **a) M and NBUE**

  $$\frac{1}{\lambda} + \min\left\{2\,(E[B] - \alpha), \frac{\rho}{2}(E[B] + \alpha)\right\}$$

  **b) Power function with parameter $c$**

  $$\frac{1}{\lambda} + \frac{(c+1)^2}{c(c+2)}E[B] - \frac{c+1}{c+2}$$

  **c) Uniform in $[0, 1]$ $(c = 1)$**

  $$\frac{1}{\lambda} + \frac{4}{3}E[B] - \frac{2}{3}$$

Finally, the ratio of the difference between the upper and the lower bound over the real value, for the exponential and power function service distributions were computed taking $\alpha = 0.5$ and $\lambda = 2, 10, 100$, and the results are in Table 3.

Table 3. Ratio of the difference between the upper and the lower bound over the real value- $\alpha = 0.5$

| Service Time Distribution | $\lambda = 2$ | $\lambda = 10$ | $\lambda = 20$ | $\lambda = 100$ |
|---|---|---|---|---|
| **Exponential** | 0.024818024 | 0.62565866 | 0.87899084 | 0.87295261 |
| **Power function with parameter $c$; $\alpha = \frac{c}{c+1}$** | 0.018536302 | 0.25071787 | 0.28865152 | 0.32992972 |

The best results (that is: the lowest) happen for the power service distribution and for the lowest traffic intensities.

## 4. Conclusions

It was already emphasized the interest of the $M|G|\infty$ system age or excess of the busy cycle, in the management of that queue, particularly of the availability of the servers.

Then this search was oriented to look for the properties of this parameter. Of course, important are the exact formulae to compute it for the various service time distributions, but some of it result quite complicate, involving infinite sums, making its applicability problematic.

So, the importance of the lower and upper bounds that it was possible to compute, mathematically much simpler, namely for service time distributions important in reliability theory such as: Exponential, NBUE, NWUE, DFR and IMRL.

**Acknowledgements**

This work was financially supported by FCT through the Strategic Project PEst-OE/EGE/UI0315/2011.

**References**

Andrade, M. (2010). A note on foundations of probability. *Journal of Mathematics and Technology, 1*(1), 96-98.

Brown, M. (1981). Further monotonicity properties for specialized renewal processes. *Annals of Probability, 9*(5), 891-895. http://dx.doi.org/10.1214/aop/1176994317

Cox, D. R. (1962). *Renewal Theory*. Methuen, London.

Ferreira, M. A. M. (1996). Valor médio da idade ou dos excesso do ciclo de ocupação na fila de espera $M|G|\infty$. *A Estatística a Decifrar o Mundo - Actas do IV Congresso Anual da Sociedade Portuguesa de Estatística, 9-15*, 231-237.

Ferreira, M. A. M., & Andrade, M. (2009). $M|G|\infty$ queue system parameters for a particular collection of service time distributions. *AJMCSR-African Journal of Mathematics and Computer Science Research, 2*(7), 138-141.

Ferreira, M. A. M., & Andrade, M. (2009a). The ties between the $M|G|\infty$ queue system transient behavior and the busy period. *International Journal of Academic Research, 1*(1), 84-92.

Ferreira, M. A. M., & Andrade, M. (2011). Fundaments of theory of queues. *International Journal of Academic Research, 3*(1, Part II), 427-429.

Ferreira, M. A. M., & Andrade, M. (2012). Busy period and busy cycle distributions and parameters for a particular $M|G|\infty$ queue system. *American Journal of Mathematics and Statistics, 2*(2), 10-15.

Ferreira, M. A. M., & Andrade, M. (2012a). Queue networks with more general arrival rates. *International Journal of Academic Research, 4*(1, PART A), 5-11.

Ferreira, M. A. M., Andrade, M., Filipe, J. A., & Coelho, M. P. (2011). Statistical queuing theory with some applications. *International Journal of Latest Trends in Finance and Economic Sciences, 1*(4), 190-195.

Hokstad, P. (1979). On the relationship of the transient behavior of a general queueing model to its idle and busy period distributions. *Mathematische Operationsforschung und Statistik. Series Optimization, 10*(3), 421-429.

Ross, S. (1983). *Stochastic Processes*. Wiley. New York.

Sathe, Y. S. (1985). Improved bounds for the variance of the busy period of the $M|G|\infty$ queue. *Advances in Applied Probability, 17*, 913-914. http://dx.doi.org/10.2307/1427096

Takács, L. (1962). *An Introduction to queueing theory.* New York: Oxford University Press.